
\baselineskip=14pt
\parskip=10pt

\magnification=\magstephalf

\def\1{{\overline{1}}}
\def\2{{\overline{2}}}
\parindent=0pt
\overfullrule=0in

\def\frac#1#2{{#1 \over #2}}

\centerline
{
 \bf On the Average Maximal Number of Balls in a  Bin Resulting from 
}
\centerline
{
\bf Throwing $r$ Balls into $n$ Bins $T$ times 
}
\rm
\bigskip
\centerline
{\it By  Amir BEHROUZI-FAR  and Doron ZEILBERGER}

{\bf Abstract.} We use the {\it holonomic ansatz} to estimate the asymptotic behavior, in $T$, of the 
average maximal number of balls in a bin that
is obtained when one throws uniformly at random (without replacement) $r$ balls into $n$ bins, $T$ times.
Our approach works, in principle, for any fixed $n$ and $r$. We were able to do the cases $(n,r)$ = $(2,1),(3,1),(4,1), (4,2)$,
but things get too complicated for larger values of $n$ and $r$.
We are pledging a \$150 donation to the OEIS for an explicit expression, 
(in terms of $n$, $r$, and $\pi$) for the constant $C_{n,r}$ such that 
that average equals $\frac{n}{r}\,T+C_{n,r} \sqrt{T}+O(1/\sqrt{T})$.

{\bf Added May 23, 2019}: Marcus Michelen came as close as possible to meeting our challenge. 
See his interesting note  {\it A Short Note on the Average Maximal Number of Balls in a Bin}, \hfill\break
{\tt https://arxiv.org/abs/1905.08933} (submitted May 22, 2019). Before that Brendan McKay
sketched a version of this argument, see \hfill\break
{\tt http://sites.math.rutgers.edu/\~{}zeilberg/mamarim/mamarimhtml/binsFeedback.html}\quad. \hfill\break
A donation of \$150 to the OEIS, in honor of Marcus Michelen and Brendan McKay, has been made.

{\bf Throwing $r$ balls into $n$ bins $T$ times}

The following problem arises in computer systems where jobs get to execute, either completely or partially, at different servers. 
Specifically, the load distribution across servers can be studied using balls and bins problem.

Suppose that you have  $n$ bins, and in each round,  you throw $r$ balls such that each ball lands in a different bin.
We assume that each of the ${{n} \choose {r}}$ possibilities are equally likely.

After  $T$ rounds, you look at the maximum occupancy amongst the $n$ bins. Let $U(n,r;T)$ be that random variable.
The smallest value it can take is, of course, $\frac{r}{n}\,T$, where all the $Tr$ balls are equally shared by all the $n$ bins,
while the largest value it can take is $T$, where one lucky  bin gets a ball at each of the $T$ rounds.

Let $A(n,r;T)$ be the expectation of $U(n,r;T)$ minus $\frac{r}{n} \, T$. Our questions are:

$\bullet$ For specified $n$ and $r$, can we compute $A(n,r;T)$ {\bf exactly} for very large values of $T$?

$\bullet$ For specified $n$ and $r$, can we determine the asymptotic behavior (in $T$), of $A(n,r;T)$?

The answer to the first question is ``{\bf of course}!'', but
only {\it in principle}. If $P(n,r,T;a_1, \dots, a_n)$ is the probability that
after $T$ rounds the occupancy of the $n$ bins are $a_1, \dots, a_n$ respectively, then
$$
P(n,r,T;a_1, \dots, a_n)=
{{{n} \choose {r}}}^{-T} \, \cdot \,
Coefficient\,\,Of\,\, x_1^{a_1} \cdots x_n^{a_n} \quad in \quad  e_{r}(x_1, \dots, x_n)^T \quad,
$$

where $e_r(x_1, \dots, x_n)$ is the {\it elementary symmetric function} of degree $r$ in $x_1, \dots, x_n$, that may be defined by
$$
e_r(x_1, \dots, x_n) = \sum_{1\leq i_1 < i_2 < \dots <i_r \leq n} x_{i_1} \cdots x_{i_r} \quad .
$$

Since we can compute $P(n,r,T;a_1, \dots, a_n)$, we can also compute
$$
Pr(U(n,r,T) \leq m) \, = \, \sum_{a_1, \dots, a_n \leq m} P(n,r,T; a_1, \dots, a_n) \quad ,
$$
from which we can compute the probability mass function, and from that, in turn, the expectation.
Unfortunately, {\it in practice}, this will not work for very large $T$, even for small $n$ and $r$.

However, using {\it Wilf-Zeilberger algorithmic proof theory} [WZ] we know {\bf for sure}, that for any specific $n$ and $r$,
our sequence of interest, $A(n,r;T), 1\leq T< \infty$ satisfies a {\it linear recurrence equation with  polynomial coefficients}
with respect to $T$. While there exist algorithms to find these recurrences from first principle, they are very slow. It is much more
efficient to compute the first few terms (using the above brute-force method), and then {\it guess} a linear
recurrence, that we can justify rigorously {\it a posteriori}, since we know that such a recurrence exists, and
we can  theoretically bound its order. We found such recurrences for the cases $(n,r)=(2,1),(3,1),(4,1),(4,2)$.

Here there are:

{\bf Proposition 1}: $A(T)=A(2,1;T)$ satisfies the following linear recurrence

{\tt A(T) = 1/(T-1)*A(T-1)+A(T-2)}, with initial conditions $A(1)=1/2,A(2)=1/2$.

In fact, it is very easy to get a {\bf closed-form} expression for $A(2T-1)$, even without WZ theory.
Since
$$
A(2T-1)\,=\, -T\,+\,\frac{1}{2} \, + \, \frac{1}{2^{2T-2}} \, \sum_{k=T}^{2T-1} \, {{2T-1} \choose {T}} T \quad ,
$$
it is readily seen that $A(2T-1)$ is given explicitly by the following  (humanly-generated!) proposition.

{\bf Proposition 1'} 
$$
A(2T-1) \,= \, \frac{(2T-1)!}{2^{2T-1} (T-1)!^2} \quad,
$$
that implies, via Stirling's formula, that
$$
A(T)=\frac{1}{\sqrt{2\pi}} \, \sqrt{T} +O(\frac{1}{\sqrt{T}}) \quad .
$$

{\bf Proposition 2}: $A(T)=A(3,1;T)$ satisfies the following fifth-order linear recurrence equation

{\tt
A(T) = 2/3*(7*T**4-32*T**3+36*T**2+19*T-36)/(7*T**2-25*T+20)/(T-1)**2*A(T-1)
+1/3*(7*T**5-53*T**4+131*T**3-55*T**2-206*T+200)/(T-2)/(7*T**2-25*T+20)/(T-1)**2*A(T-2)
+1/3*(21*T**4-180*T**3+519*T**2-544*T+160)/(7*T**2-25*T+20)/(T-1)**2*A(T-3)
}
{\tt-2/3*(T-3)*(7*T**3-32*T**2+29*T+2)/(7*T**2-25*T+20)/(T-1)**2*A(T-4)
-1/3*(T-3)*(T-4)*(7*T**2-11*T+2)/(7*T**2-25*T+20)/(T-1)**2*A(T-5)
},   \hfill\break
with initial conditions $A(1)=2/3$, $A(2)=2/3$, $A(3)=8/9$, $A(4)=28/27$, $A(5)=10/9$.

{\bf Proposition 3}: $A(T)=A(4,1;T)$ satisfies a ninth-order linear recurrence equation that can be found here:
{\tt http://sites.math.rutgers.edu/\~{}zeilberg/tokhniot/oBINnr3.txt} \quad .

{\bf Proposition 4}: $A(T)=A(4,2;T)$  satisfies an eighth-order linear recurrence equation that can be found here:
{\tt http://sites.math.rutgers.edu/\~{}zeilberg/tokhniot/oBINnr4.txt} \quad .

It also follows (at least in the above cases, but probably can be proved in general), that for each $n$ and $r$ there
exists a constant $C_{n,r}$ such that, asymptotically
$$
A(n,r;T) = C_{n,r} \sqrt{T} +O(\frac{1}{\sqrt{T}}) \quad .
$$
By using the recurrences to compute many values, we were able to derive the following estimates
$$
C_{2,1}=0.3989\dots \quad, \quad C_{3,1}= 0.489\dots \quad , \quad C_{4,1}= 0.516 \dots \quad, \quad C_{4,2}= 0.59430 \dots \quad .
$$
(Note that the value for $C_{2,1}$ agrees with the human-generated {\it exact value} given by Proposition $1'$, $(2 \pi)^{-1/2}=0.398942280401\dots$).

Using the algorithm in [Z], it is possible to get higher-order asymptotics, that usually enable a much more
precise estimate of the leading constant of the asymptotics. To our disappointment, this does not seem
to be the case for this problem, hence we do not bother to state this more refined asymptotics.
Our computers are not large enough to handle larger values of $n$ and $r$, even though we do know that such
a recurrence exists for every choice of $n$ and $r$.

{\bf The Maple package.} The accompanying  Maple package  {\tt BINnr.txt} is available from the front of this article:\hfill\break
{\tt http://sites.math.rutgers.edu/\~{}zeilberg/mamarim/mamarimhtml/bins.html} \quad .

{\bf Estimating $C_{n,r}$ for arbitrary $n$ and $r$}

Using heuristic arguments,  Behrouzi-Far and Soljanin  obtained estimates for $C_{n,r}$, that do not quite agree with the above
values  for small $n$ and $r$, but are hopefully close to the exact values for larger $n$ and $r$. It is

$$
(r/n)\,\sqrt{\pi \, log(n)}\,\log(n/r) \quad .
$$

{\bf Conclusion.} While the results of this paper may be disappointing to people who are interested in
this problem for large $n$ and $r$, it is interesting from the point of view of {\it experimental mathematics}.
By using the starting values of the sequence of interest, $A(n,r;T)$, as `training data', that were computed
by brute force, the computer `learned' the linear recurrence, by fitting it into an {\it ansatz} (see [KP], chapter 7)
and this enabled it to {\it predict} (exactly!) the values for very large $T$, that would be impossible to
compute, in practice, by the naive method.

{\bf Encore: An asymptotic challenge}

Our quantity of interest, $A(n,r;T)$, is expressible as a multi-variable contour integral,
hence it may be  possible to find the asymptotic constant, $C_{n,r}$, as a closed-form expression
in terms of $n$ and $r$ and $\pi$. 
It is also expressible, via the multinomial theorem, as a huge multi-sum featuring multinomial coefficients.
The special case $r=1$ `only' has $n-1$ $\sum$-signs, and may be not too difficult, but
the general $r$ seems harder. The authors of the modern classic [PW], or one of their disciples, may rise to the challenge of finding that expression.

We are  pledging a \$50 donation to the OEIS in honor of the first person to find a closed-form expression,
in terms of $n$ and $\pi$, for the case $r=1$, i.e. for $C_{n,1}$, and a further \$100 donation for the general
case, i.e. for $C_{n,r}$.

{\bf Acknowledgment.} We thank  Emina Soljanin for introducing us to this interesting problem.

{\bf References}

[BS] Amir Behrouzi-Far and Emina Soljanin, {\it  On the load distribution in systems with bimodal job sizes}, preprint.

[KP] Manuel Kauers  and Peter Paule, {\it ``The Concrete Tetrahedron''}, Springer, 2011.

[PW] Robin Pemantle and Mark C. Wilson, {\it ``Analytic Combinatorics in Several Variables''},
Cambridge University Press, 2013.

[WZ] Herbert S. Wilf and Doron Zeilberger,{\it An algorithmic proof theory for hypergeometric (ordinary and "q") multisum/integral identities},
 Invent. Math. {\bf 108}(1992), 575-633. Available from \hfill\break
{\tt http://sites.math.rutgers.edu/\~{}zeilberg/mamarim/mamarimhtml/multiwz.html} \quad .

[Z] Doron Zeilberger,
{\it AsyRec: A Maple package for Computing the Asymptotics of Solutions of Linear Recurrence Equations with Polynomial Coefficients},
The Personal Journal of Shalosh B. Ekhad and Doron Zeilberger,  April 4, 2008.  \hfill\break
{\tt http://sites.math.rutgers.edu/\~{}zeilberg/mamarim/mamarimhtml/asy.html} \quad .

\bigskip

\hrule
\bigskip
 Amir Behrouzi-Far, Department  of  Electrical  and  Computer   Engineering,   Rutgers   University,   Piscataway,   NJ   08854,   USA;
amir dot behrouzifar at rutgers dot edu
\bigskip
Doron Zeilberger, Department of Mathematics, Rutgers University (New Brunswick), Hill Center-Busch Campus, 110 Frelinghuysen
Rd., Piscataway, NJ 08854-8019, USA. \hfill \break
DoronZeil at gmail dot com  \quad ;  \quad {\tt http://sites.math.rutgers.edu/\~{}zeilberg/} \quad .
\bigskip
\hrule
\bigskip
Exclusively published in The Personal Journal of Shalosh B. Ekhad and Doron Zeilberger  \hfill \break
{ \tt http://sites.math.rutgers.edu/\~{}zeilberg/pj.html} \quad  and {\tt arxiv.org} \quad . 
\bigskip
\hrule
\bigskip

{\bf First  Written:  May 20, 2019. This version: May 23, 2019.} 

\end